\documentclass[11pt]{amsart}
\usepackage{amscd,amssymb,amsmath}
\usepackage{graphicx,epsfig}
\newtheorem{theorem}{Theorem}[section]
\newtheorem{lemma}[theorem]{Lemma}
\newtheorem{proposition}[theorem]{Proposition}

\theoremstyle{definition}

\theoremstyle{remark}
\newtheorem{remark}{Remark}[section]

\numberwithin{equation}{section}
\begin{document}
\title{Equivariant path fields on topological manifolds}
\author{Luc\'\i lia Borsari}
\address{Departamento de Matem\'atica - Instituto de Matem\'atica e Estat\'istica - Universidade de S\~ao Paulo, Rua do Mat\~ao, 1010 - CEP 05508-090, S\~ao Paulo - SP, Brasil}
\email{lucilia@ime.usp.br}

\author{Fernanda Cardona}
\address{Departamento de Matem\'atica - Instituto de Matem\'atica e Estat\'istica - Universidade de S\~ao Paulo, Rua do Mat\~ao, 1010 - CEP 05508-090, S\~ao Paulo - SP, Brasil}
\email {cardona@ime.usp.br}

\author{Peter Wong}
\address{Departament of Mathematics - Bates College - Maine, USA}
\email {pwong@bates.edu}
\begin{abstract}
A classical theorem of H. Hopf asserts that a closed connected smooth manifold admits a nowhere vanishing vector field if and only if its Euler characteristic is zero. R. Brown generalized Hopf's result to topological manifolds, replacing vector fields with path fields. In this note, we give an equivariant analog of Brown's theorem for locally smooth $G$-manifolds where $G$ is a finite group.
\end{abstract}
\date{\today}
\thanks{The third author acknowledges supported by a grant from the National Science Foundation.}

\newcommand{\bea} {\begin{eqnarray*}}
\newcommand{\beq} {\begin{equation}}
\newcommand{\bey} {\begin{eqnarray}}
\newcommand{\eea} {\end{eqnarray*}}
\newcommand{\eeq} {\end{equation}}
\newcommand{\eey} {\end{eqnarray}}

\newcommand{\ovl}{\overline}
\newcommand{\vv}{\vspace{4mm}}
\newcommand{\lra}{\longrightarrow}
\newcommand{\bg}{\bigskip}
\newcommand{\md}{\medskip}
\newcommand{\sm}{\smallskip}
\newcommand{\oo}{\hbox{\itshape\bfseries o}}



\maketitle

\section{Introduction}

Let $M$ be a closed connected orientable smooth manifold. A classical theorem of H. Hopf \cite{hopf} states that $M$ admits a non-singular vector field if and only if the Euler characteristic, $\chi(M)$, of $M$ is zero. R. Brown \cite{brown} gave a generalization of Hopf's theorem for topological manifolds, by replacing vector fields with path fields, a concept first introduced by J. Nash \cite{nash}. In \cite{brown}, R. Brown showed that a compact topological manifold $M$ admits a non-singular path field if and only if $\chi(M)=0$. Subsequently, R. Brown and E. Fadell \cite{brown-fadell} extended \cite{brown} to topological manifolds with boundary.  It was shown by E. Fadell \cite{fadell3} that any Wecken complex of zero Euler characteristic admits a non-singular simple path field. R. Stern \cite{stern} showed the same result for topological manifolds of dimension different from four.

The existence of a path field allows one to show the so-called {\it Complete Invariance Property} (CIP) (see \cite{jiang-schirmer} and \cite{schirmer}). Recall that a topological space $M$ is said to have the CIP if for any non-empty closed subset $A\subset M$, there exists a map $f:M\to M$ such that $A=Fix\,f:=\{x\in M\,|\,f(x)=x\}$. Similarly, $M$ possesses the CIP with respect to deformation (denoted by CIPD) if $f$ is homotopic to the identity $1_M$. The non-singular path field problem is equivalent to the fixed point free deformation problem. That is, $M$ admits a non-singular path field if and only if $1_M$ is homotopic to a fixed point free map.

In \cite{komiya1}, \cite{komiya2}, and \cite{wil}, equivariant vector fields on compact smooth $G$-manifolds were studied. In particular, an equivariant analog of Hopf's theorem was proved in \cite{komiya1}. Furthermore, an equivariant analog of what was done for path fields on Wecken complexes in \cite{fadell3}, was given in \cite{wong3} and necessary and sufficient conditions for equivariant CIPD were given for smooth $G$-manifolds (see also 
\cite{azad-srivastava} for a certain type of equivariant CIP). Similar to the non-equivariant case, the equivariant non-singular path field problem is closely related to finding an equivariant fixed point free deformation. It turns out that the existence of such a fixed point free map requires more than merely the existence of non-equivariant fixed point free deformation on the fixed point sets $M^H$ for each isotropy type $(H)$ (see \cite{ferrario}).

The objective of this paper is to prove an equivariant analog of Brown's theorem \cite{brown} for topological manifolds with locally smooth action of a finite group $G$. Moreover, we extend the necessary and sufficient conditions for $G$-CIPD found in \cite{wong4} to this category of $G$-manifolds.

We would like to thank D.L. Gon\c calves and G. Peschke for very helpful conversations and suggestions.

Throughout $G$ will always be a finite group acting on a compact space $M$ where the action is locally smooth. For the definition and basic properties of locally smooth actions, we refer the reader to \cite{bredon}.

\section{Equivariant Euler characteristic and $G$-path fields}

In this section, we establish the necessary definitions of path fields and Euler characteristic in the equivariant category.

Equivariant path fields were defined and studied in \cite{wong3} and \cite{wong4}. For our purposes, we think of $G$-path fields as sections of certain $G$-fibrations.

First, given a $G$-map $p:E\to B$, we say that $p$ has the $G$-{\it Covering Homotopy Property} ($G$-CHP) if for all $G$-space $X$ the following commutative diagram has a solution $F:X\times [0,1] \to E$ where all maps are $G$-equivariant.

\begin{equation*}
\begin{CD}
X \times \{0\}@>f>> E \\
@V{\text{incl.}}VV  @VVpV\\
X \times [0,1] @>H>> B
\end{CD}
\end{equation*}
A {\it $G$-fibration} is simply a $G$-map $p:E\to B$ satisfying the $G$-CHP for all $G$-spaces. 

Given a $G$- fibration $p:E \to B$, we consider $\Omega _p = \lbrace (e,w)\in E
  \times B^I | p(e) = w(0)\rbrace$. Then $\Omega_p$ is a $G$-invariant
  subspace of $ E\times B^I$. Let $\widetilde p :E^I \to \Omega_p$ be the $G$-map
  defined by  $\widetilde p(\tau) = (\tau(0), p(\tau))$. Consider the equivariant maps
  $F: \Omega_p \times I \to B$ defined by $F(e,w,t) = w(t)$, and $f:\Omega_p \to E$ by $f(e,w) = e$. Since $p$ is a $G$-fibration, $F$ can be lifted to  a $G$-map $\widetilde F:\Omega_p \times I \to B$ which extends $f$. Then $\lambda:\Omega_p \to E^I$, defined by $\lambda(e,w)(t) = \widetilde F(e,w,t)$, is an equivariant lifting function for $p$, that is, $\widetilde p \circ \lambda$ is the identity on $\Omega_p$.
  
  A $G$-fibration is called regular if it admits a regular $G$-lifting function, meaning, a $G$-lifting function satisfying $\lambda(e,\overline{p(e)}) = e$, for all $e\in E$, where $\overline{p(e)}$ denotes the constant path at $p(e)$. In \cite{hurewicz}, W. Hurewicz shows that every fibration over a metric space is regular. The same proof can be adapted to the equivariant case, provided the metric $d$ is assumed to be $G$-invariant, that is, $d(gx,gy) = d(x,y)$, for all $g \in G$ and $x,y \in B$. 

\begin{lemma}\label{lemma1}
Let $p:E \to B$ be a regular $G$-fibration over a $G$-manifold $B$. Let $(X,A)$ be a $G$-ANR pair and suppose that there are equivariant maps $f: X\times {0} \cup A\times I \to E$ and $h: X \times I \to M$ such that $p \circ f = h|_{(X\times {0} \cup A\times I)}$. Then, there exists a $G$-map $\widetilde f : X \times I \to E$ which extends $f$ and such that $p \circ \widetilde f = h$.
\end{lemma}
\begin{proof}
This lemma is an equivariant version of Theorem 2.4 of \cite{allaud-fadell}. The proof of this theorem in the non equivariant context is very constructive and it is possible to verify that, in all steps, we do obtain equivariant maps, as long as we start with the appropriate equivariant setting and make use of Corollary 2.3 of \cite{wil}.
\end{proof}

Given a compact topological manifold $M$, the {\it Nash path space} $ T_M$ of $M$ consists of $\overline { T_M}=\{\text{all~constant~paths}\}$ and the set $ T^0_M$ of all paths $\alpha$ on $M$ such that for $0\le t\le 1$, $\alpha(t)=\alpha(0)$ iff $t=0$. Consider the map $q: T_M \to M$ given by $q(\alpha)=\alpha (0)$. With the compact-open topology on $ T_M$, the triple $( T_M, q_M,M)$ is a Hurewicz fibration and the sections of $q$ are called {\it path fields} on $M$. A path field is 
{\it non-singular} if it is a section in $(T^0_M, q_M|_{ T^0_M},M)$. A path field $\sigma$ is 
{\it simple} if for any $x\in M$, $\sigma(x)$ is a simple path.

If $G$ acts on $M$, then $G$ acts on $ T_M$ via
$g*\alpha(t)=g\alpha(t)$. Since $q: T_M \to M$ is a fibration, it is
straightforward to see that it is indeed a $G$-fibration where the $G$-action
on $[0,1]$ is trivial. Thus, we define a {\it $G$-path field} to be a $G$-section
$s:M\to  T_M$ of $q$ so that $q\circ s=1_M$. Moreover, the
subfibration $q_M^0: T^0_M \to M$ is also a $G$-subfibration. The
notions of non-singular and of simple $G$-path fields are defined in the
obvious fashion.

Given a compact topological manifold $M$, the classical Euler characteristic
of $M$ is an integer and it coincides with the fixed point index of the
identity map $1_M$. When a finite group $G$ acts on $M$, the appropriate
equivariant Euler characteristic takes the components of the various fixed
point sets $M_H$, $H\leq G$, into account. 

We write $|\chi| (M_H)=\sum_{C} |\chi (C)|$, where
$C$ ranges over the connected components of $M_H=\{x\in M|G_x=H\}$. Here, $G_x$
denotes the isotropy subgroup of $x$. Since $M$ is compact, each $M^H=\{x\in M|hx=x,\ \forall h\in M\}$ is also compact so that $M_H$ has only a finite number of components.

\section{Singularities of $G$-path fields}

In this section, we prove our main results following the approach of \cite{brown}. Since we work in the $G$-manifolds category, many of the techniques employed in \cite{brown} must be modified for the equivariant setting, first of which is the following relative equivariant domination theorem for compact $G$-ANRs.
 
\begin{theorem}[Relative Equivariant Domination Theorem]\label{dominate}
Let $M$ be an $n$-dimensional
$G-$manifold and $A$ be an invariant compact submanifold of dimension $k$. We can find a $G-$complex $K$ of dimension $n$, an invariant subcomplex of dimension $k$ and
 equivariant maps $\varphi\colon K\lra M$ and  $\psi\colon M\lra K$,
so that $\psi$ is barycentric, $\varphi|_L\colon L\lra A$, $\psi|_A\colon A\lra L$, $ \varphi\circ\psi{\cong}_{G}\ \hbox{id}_M$ and $\varphi|_L\circ\psi|_A{\cong}_{G}\ \hbox{id}_A$
\end{theorem}
\begin{proof}
According to \cite[Theorem 1]{antonian}, we can equivariantly embed $M$ as a closed $G$-neighborhood retract of a convex $G$-set in a Banach $G$-space $A(M)$ in which $G$ acts isometrically. Now we follow the proof of the $G$-domination theorem (Proposition 2.3) of \cite{kwasik}. Let $r:{\mathcal O}\to M$ be the $G$-retraction of some $G$-invariant neighborhood ${\mathcal O}$. It suffices to show that $\mathcal O$ can be $G$-dominated by a finite $G$-complex $K$. Let $\{W_{\alpha}\}$ be a covering of $\mathcal O$ by convex subsets which are open in $\mathcal O$. Since $M$ is compact, we can find a finite open covering $\{\mathcal O_{\alpha}\}$ of $\mathcal O$ such that the convex hull of a finite union of the $\mathcal O_{\beta}$ is contained in $M$. Then there is a finite open pointed $G$-covering $\mathcal V=\{V_{\gamma},v_{\gamma}\}$, a refinement of $\{\mathcal O_{\beta}\}$ such that $v_{\gamma}\in A$ if $V_{\gamma}\cap A\ne \emptyset$. Let $K=|N(\mathcal V)|$ be the nerve of $\mathcal V$ with the canonical $G$-action. 

For any $x\in M$, we let
$$
\nu(x)=\sum_i d(x,M-V_i)
$$
where $d$ denotes the metric on $M$ which is $G$-invariant since $G$ acts isometrically. Let $\{v_i\}$ be the vertices of $|N(\mathcal V)|$. Now $\{\sum_i \frac{d(x,M-V_i)}{\nu(x)}\}$ is a $G$-partition of unity subordinate to $\mathcal V$. Define the $G$-map $\varphi :M \to |N(\mathcal V)|$ by 
$$
\varphi (x)=\sum_i \frac{d(x,M-V_i)}{\nu(x)}v_i.
$$
Note that $\varphi$ is a barycentric mapping. Consider the $G$-map $\psi=r\circ \eta : |N(\mathcal V)| \to M$, where $\eta$ is the map $\psi$ as in the proof of Proposition 2.3 of \cite{kwasik}. It is straightforward to check that the $G$-maps $\varphi$ and $\psi$ yield the desired $G$-domination. Note that $K$ is of dimension less than or equal to $n$ since the $\mathcal V$ is a refinement. It follows that $K$ must be of dimension $n$ otherwise $K$ has no homology in dimension $n$ whereas $\dim M=n$ and $M$ is a compact manifold of dimension $n$. Finally, we let $\mathcal V_A=\{(V_i\cap A,v_i)\}$. It follows that $\mathcal V_A$ is a $G$-covering of $A$ and the nerve $L= |N(\mathcal V_A)|$ is a subcomplex of $K$. Since $A$ is a compact manifold of dimension $k$, we conclude that $L$ is of dimension $k$ and that $L$ equivariantly dominates $A$. 
\end{proof}

\begin{remark} It has been noted by S. Antonyan in \cite{antonian} that the equivariant embedding theorem 
\cite[Theorem 6.2]{murayama} of M. Murayama is incorrect: in that the Banach space $B(M)$ of all bounded continuous functions on $M$ used in \cite{murayama} is not a Banach $G$-space and the $G$-action defined there is not continuous. Likewise, the same mistake was also committed by S. Kwasik in \cite{kwasik}. Nevertheless, the $G$-domination theorem in both \cite{kwasik} and \cite{murayama} is stated correctly and their proofs are valid provided one replaces $B(M)$ with the linear subspace $A(M)$ of all $G$-uniform functions as in 
\cite{antonian}. We thank M. Golasi\'nski for bringing \cite{antonian} to our attention. As noted by Hanner in~\cite{hanner}, in non-equivariant settings Borsuk showed in~\cite{borsuk} that any compact ANR is dominated by a finite polyhedron. Then, in~\cite{brooks}, Brooks showed, again in the non-equivariant setting, that if an $n$-dimensional compact ANR is dominated by a complex then it is dominated by its $n$-dimensional skeleton.
\end{remark}

In order to prove the next proposition, we will need the following non-equivariant result.

\begin{lemma}\label{lemma}
Let $M$ be a $n$-manifold and $A\subset M$ a submanifold. Consider $c$ an $n$-cell in $M-A$, with $\partial c\subset A$ and let $\sigma\colon\partial c\lra T_A$ be a path field  and $\oo $ be a point in the interior of $c$. Then, there exists a path field $\sigma^\prime\colon c\lra T_M$, extending $\sigma$ with $\oo$ being its only singularity in the interior of $c$. Moreover, in case $\sigma$ has singularities then we may take $\sigma^\prime$ without singularities in the interior of $c$.
\end{lemma}

\sm
\begin{proof}
 Let $c^\prime$ be an $n$-cell contained in Int$\,c$, with $\oo$ in its interior. We will extend $\sigma$ to $c-\hbox{Int}\, c^\prime$ without creating new singularities: Let $[\oo, b_x]$ be the oriented segment through $x$, beginning at $\oo$, ending at $b_x\in\partial c$. Therefore we could write any $x\in c-\hbox{Int}\, c^\prime$ as $x=(1-t_x)\oo+t_x b_x$, where $t_x\in]0,1]$. 

\bg

\includegraphics[height=5.5cm]{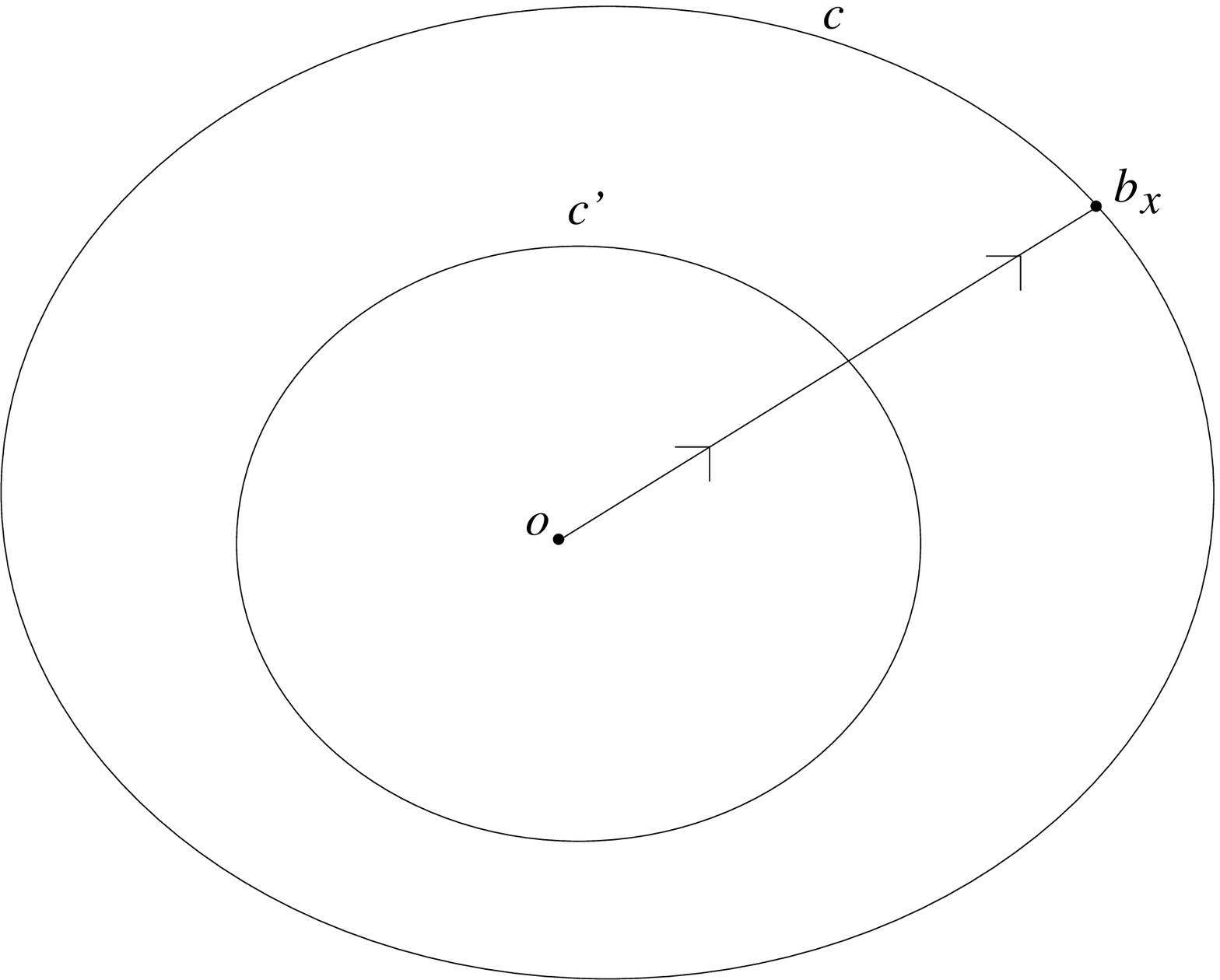}

The extended path field $\sigma^\prime$ will be defined for each $x\in c-\hbox{Int}\, c^\prime$ as follows:
$$
\sigma^\prime(x)(s)=\left\{
\begin{array}{ll}
 (1-t_x-s)\oo+(t_x+s)\,b_x,& 0\leq s\leq 1-t_x\\
&\\
\sigma(b_x)\left(\displaystyle\frac{s+t_x-1}{t_x}\right),& 1-t_x\leq s\leq 1\\
\end{array}
\right.
$$

Observe that it is well defined because $t_x>0$, for any  $x\in c-\hbox{Int}\, c^\prime$. Also, in the first equation, for any $t_x$ when $s=0$ we have $(1-t_x)\oo+t_xb_x=x$; when $s=1-t_x$ we have $(t_x+1-t_x)\,b_x=b_x$. In the second equation, when $s=1-t_x$ we have $\sigma(b_x)(\frac{1-t_x+t_x-1}{t_x})=\sigma(b_x)(0)=b_x$; when $s=1$ we have $\sigma(b_x)(\frac{1+t_x-1}{t_x})=\sigma(b_x)(1)$. Therefore $\sigma^\prime$ has no other singularities than those of
$\sigma$ (if it has any, they will be in the boundary of $c$), so $\sigma^\prime$ has no singularities in the boundary of $c^\prime$. By Lemma 1.5 of \cite{brown} it can be extended to $\hbox{Int}\, c^\prime$ having only $\oo$ as a singularity in 
$\hbox{Int}\, c^\prime$.

By an abuse of notation we will denote this extension of $\sigma^\prime$ to $\hbox{Int}\, c^\prime$ also by $\sigma^\prime$. Therefore, we constructed an extension of $\sigma$, $\sigma^\prime\colon c\lra T_M$, which has only one singularity in $\hbox{Int}\, c$ and in the boundary only the singularities that $\sigma$ had.

If $\sigma$ does have singularities in the boundary of the cell, we will eliminate the singularity of 
$\sigma^\prime$ in its interior:

Let $y\in\partial c$ be a singularity of $\sigma$ (therefore a singularity of $\sigma^\prime$); let $c_1\subset c_2\subset c$ be two cells such that $\partial c_1\cap\partial c_2=\{y\}$, $\partial c_i\cap\partial c=\{y\}$, for $i=1,2$ and $\oo\in \hbox{Int}\, c_1$ (and therefore $\oo\in \hbox{Int}\, c_2$).

Let $[b_x,y]$ be the oriented segment through $x$, beginning at $b_x\in\partial c_2$, ending at the singularity $y\in\partial c$. Therefore we could write any $x\in c_2$ as $x=(1-t_x)b_x+t_x y$, where $t_x\in[0,1]$. Also, each of these segments would determine a point $a_x \in\partial c_1$ such that $a_x=(1-\ovl t_x)b_x+\ovl t_x y$, with $\ovl t_x\in\,]0,1[$.

\vskip0.5truecm

\includegraphics[height=5.5cm]{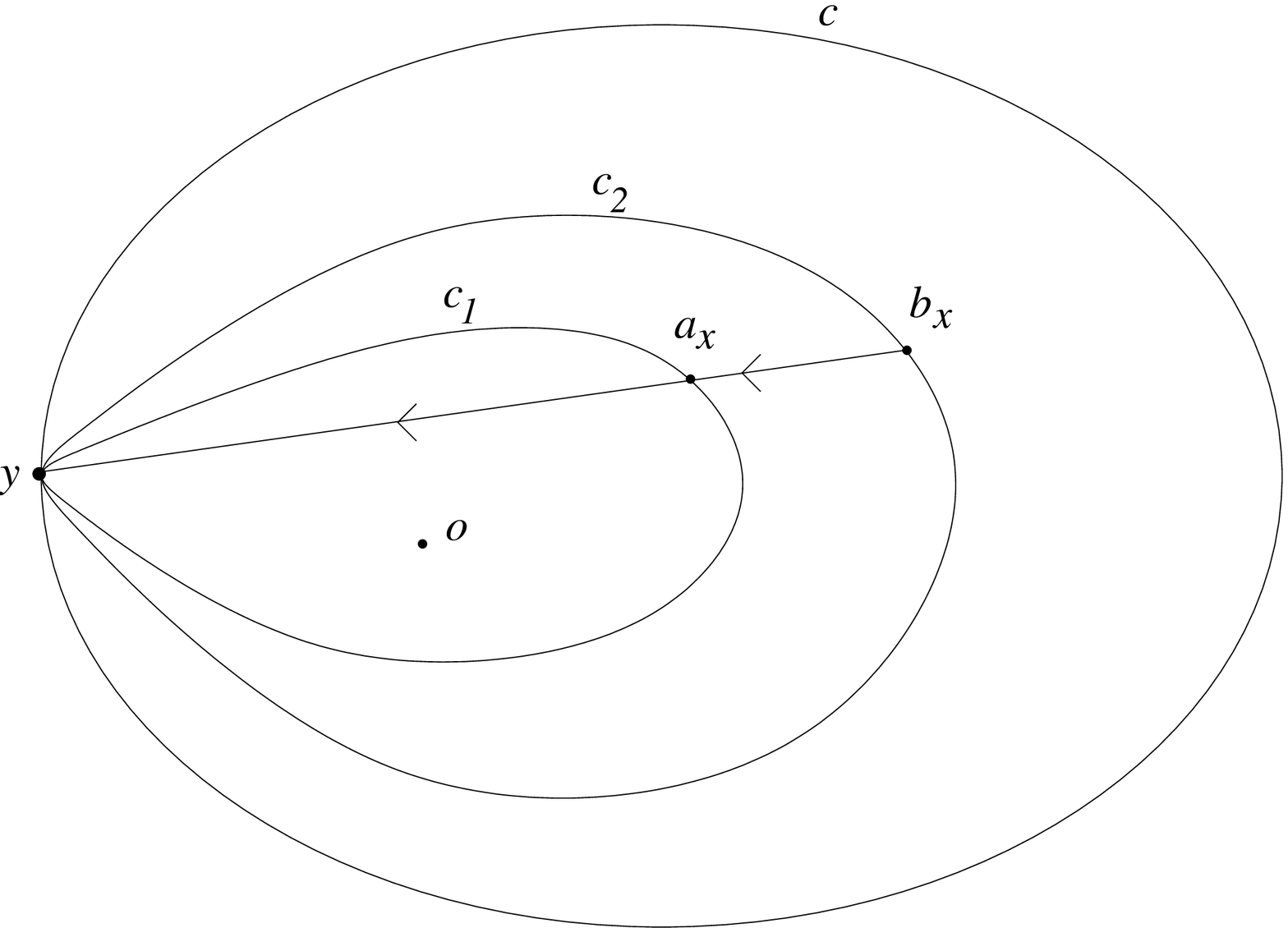}

The new path field $\ovl\sigma$ will be defined in each one of the regions represented above, as follows:

-- For $x\in\hbox{Int}\, c_1$, we have $0<t_x<1$ and 
$$ \ovl\sigma(x)(s)=\left\{\begin{array}{ll}
                            (1-t_x-s)b_x+(t_x+s)y,& 0\leq s\leq 1-t_x\\
			     &\\
			     \sigma^\prime(y)\left(\displaystyle\frac{s+t_x-1}{t_x}\right),& 1-t_x\leq s\leq 1\\
                           \end{array}
			\right.
$$

Observe that it is well defined because $t_x>0$, for any  $x\in \hbox{Int}\, c_1$. Also, in the first equation, for any $t_x$ when $s=0$ we have $(1-t_x)b_x+t_xy=x$; when $s=1-t_x$ we have $(t_x+1-t_x)\,y=y$. In the second equation, when $s=1-t_x$ we have $\sigma^\prime(y)(\frac{1-t_x+t_x-1}{t_x})=\sigma^\prime(y)(0)=y$.

\sm

-- For $x\in c_2-\hbox{Int}\, c_1$, we have $0\leq t_x\leq \ovl t_x $ and 
$$ \ovl\sigma(x)(s)=\left\{\begin{array}{ll}
                            (1-t_x-s)b_x+(t_x+s)y,& 0\leq s\leq \displaystyle\frac{t_x}{\ovl t_x}-t_x\\
			     &\\
    \sigma^\prime\left(\left(1-\displaystyle\frac{t_x}{\ovl t_x }\right)b_x+\displaystyle\frac{t_x}{\ovl t_x}y\right)
\left(\displaystyle\frac{s-(\frac{t_x}{\ovl t_x}-t_x)}{1-(\frac{t_x}{\ovl t_x }-t_x)}\right),
& \displaystyle\frac{t_x}{\ovl t_x}-t_x\leq s\leq 1\\
                           \end{array}
			\right.
$$

Observe that it is well defined because $\ovl t_x>0$, for any  $x\in c_2$. Also, in the first equation, for any $t_x$ when $s=0$ we have $(1-t_x)b_x+t_x y=x$; when $s=\frac{t_x}{\ovl t_x}-t_x$ we have $(1-t_x-(\frac{t_x}{\ovl t_x}-t_x))\,b_x+(t_x+\frac{t_x}{\ovl t_x}-t_x)\,y= (1-\frac{t_x}{\ovl t_x})\,b_x+\frac{t_x}{\ovl t_x}\,y$. In the second equation, when $s=\frac{t_x}{\ovl t_x}-t_x$ we have $\sigma^\prime((1-\frac{t_x}{\ovl t_x })b_x+\frac{t_x}{\ovl t_x}y)(0)=(1-\frac{t_x}{\ovl t_x })b_x+\frac{t_x}{\ovl t_x}y$.
\sm

-- For $x\in c-\hbox{Int}\, c_2$, 
$$ \ovl\sigma(x)=\sigma^\prime(x)\ .$$

\sm

Notice that if $x\in\partial c_1$ then $t_x=\ovl t_x$ and if $x\in\partial c_2$ then $t_x=0$ and therefore $\ovl\sigma$ is well-defined and continuous in $\partial c_1$ and $\partial c_2$, the boundaries of $c_1 $ and $c_2$. A simple verification will show that $\ovl\sigma$ has no singularities in $\hbox{Int}\, c$ and the fact that $\sigma$ is a path field in $A$ guarantees that $\ovl\sigma$ is in fact a path field.

\end{proof}

\begin{proposition}\label{relative free case}
Let $M$ be a locally smooth $G-$manifold, $\dim M=n$,
$A\subset M$ an invariant submanifold so that $G$ acts freely on
$M-A$. Given an equivariant section $\sigma_A\colon A\lra T_A$, with a finite number of singular orbits, it is possible to extend $\sigma_A$ to an equivariant section $\sigma\colon M\lra T_M$ in such a way the closure of each component of $M-A$ intersects at most one singular orbit of $\sigma$.
\end{proposition}
\begin{proof}
Consider the following diagram:

\bg

\includegraphics[height=8cm]{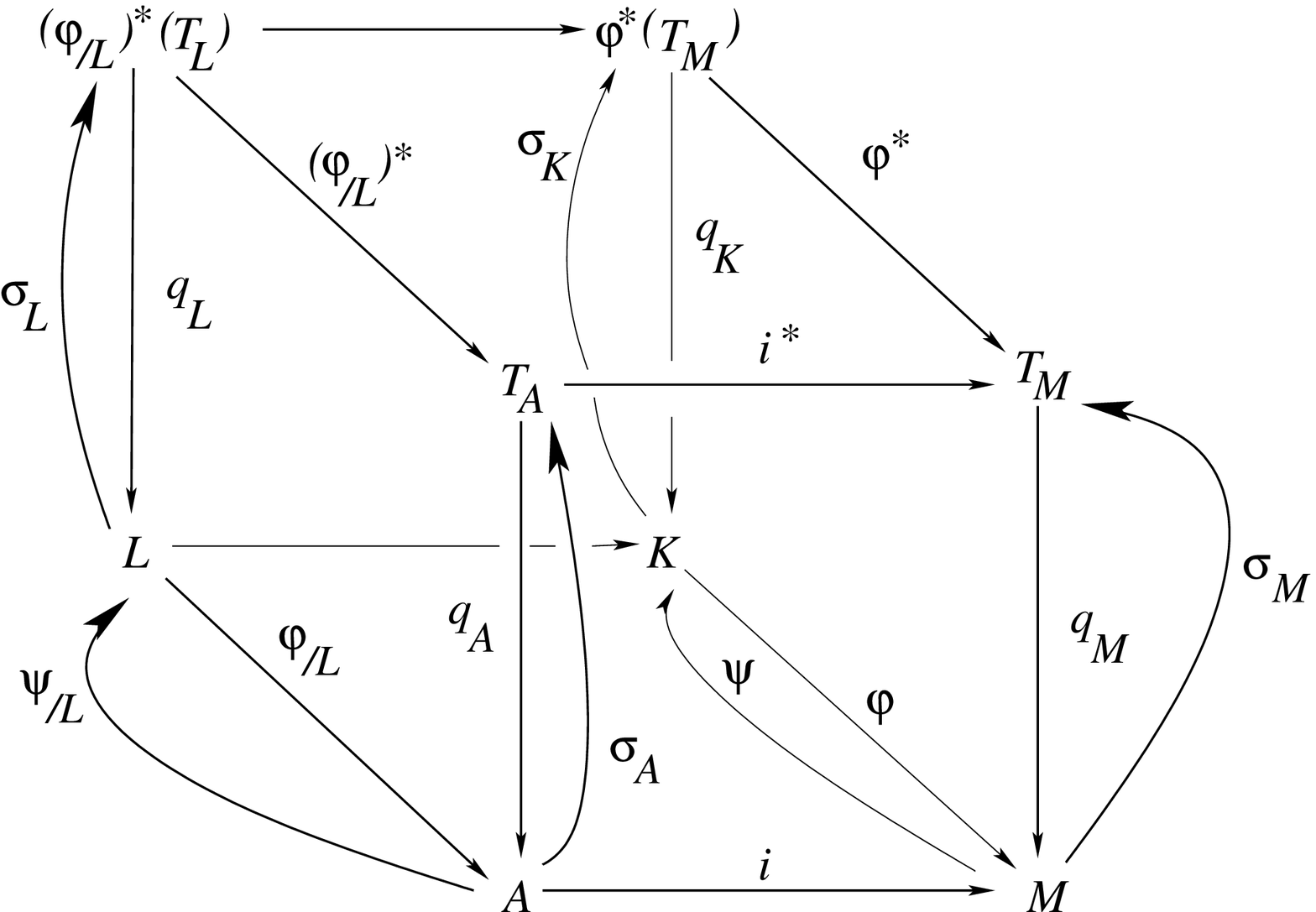}

\bg
Here, $K$ and $L$ are as in Theorem~\ref{dominate}  and $q_L\colon (\varphi|_L)^* (T_L)\longrightarrow L$ and $q_K\colon \varphi^* (T_M)\longrightarrow K$ denote the pullbacks of $q_A$ and $q_M$ by $\varphi$, respectively.

Now, starting with $\sigma_A$, a $G$-section of $q_A:T_A \to A$, we  define a $G$-section
$\sigma_L: L \to (\varphi|_L)^*(T_A)$, by $\sigma_L(y) = (y, \sigma_A(\varphi(y))$. A similar procedure as the one indicated in Lemma 1.6 in \cite{brown} can be used to extend $\sigma_L$  to a $G$-section $\sigma_K$, having only a  finite number of singular orbits in
$K-L$. In order to extend $\sigma_L$ to an $m$-simplex $\delta$ of $K-L$, we use Lemma~\ref{lemma}  and extend it to $g\delta$ in the usual equivariant
way. Define $\sigma'_M :M \to T_M$ by $\sigma'_M(x) =
\varphi(\sigma_K(\psi(x))$. Then $\sigma'_M$ is an equivariant map, but it is
not a section. In fact, it is a homotopy section since
$q_M(\sigma_M(x))=\varphi\circ\psi(x)$. Consider $h:M\times I \to M$, the $G$-homotopy
between $\varphi \circ \psi$ and the identity on $M$, and $f:A\times I \cup M\times{0}\to T_M$ given by the $G$-homotopy between $\sigma_A$ and $(\sigma'_M)|_A$ on $A\times I$ and by $\sigma'_M$ in $M\times{0}$. Since $M$ admits an invariant metric, we may apply
Lemma~\ref{lemma1} to obtain an equivariant lifting $f'$ of $h$ extending $f$. Define
$\sigma_M :M \to T_M$ by $\sigma_M(x) = f'(x,1)$.  Then $\sigma_M$ is a
$G$-section on $M$ extending $\sigma_A$.

The first step is to change  $\sigma_M$ to reduce the singular set in $M-A$ to a finite
one. In order to do this, consider first $\{Gx_1,Gx_2,...,Gx_r\}$ the set of singular
orbits of $\sigma_K$ in $K-L$. The set of singular orbits of  $\sigma_M$, which are not in $A$, lies in
the pre-image of $\{Gx_1,Gx_2,...,Gx_r\}$ under $\psi$.  Since $\psi$ is
equivariant, this set is
$\{G\psi^{-1}(x_1),...,G\psi^{-1}(x_r)\}$. Since $G$ acts freely in $M-A$, for each
$i$, the sets $g\psi^{-1}(x_i)$, $g$ in $G$, are disjoint. Following the proof
of Theorem 1.10 of \cite{brown}, we can assume that for a connected component $C$ of $M-A$, $\psi^{-1}(x_1)\cap C$ is contained in the interior of $c$, a closed topological $n$-cell (see the figure below).

\includegraphics[height=9cm]{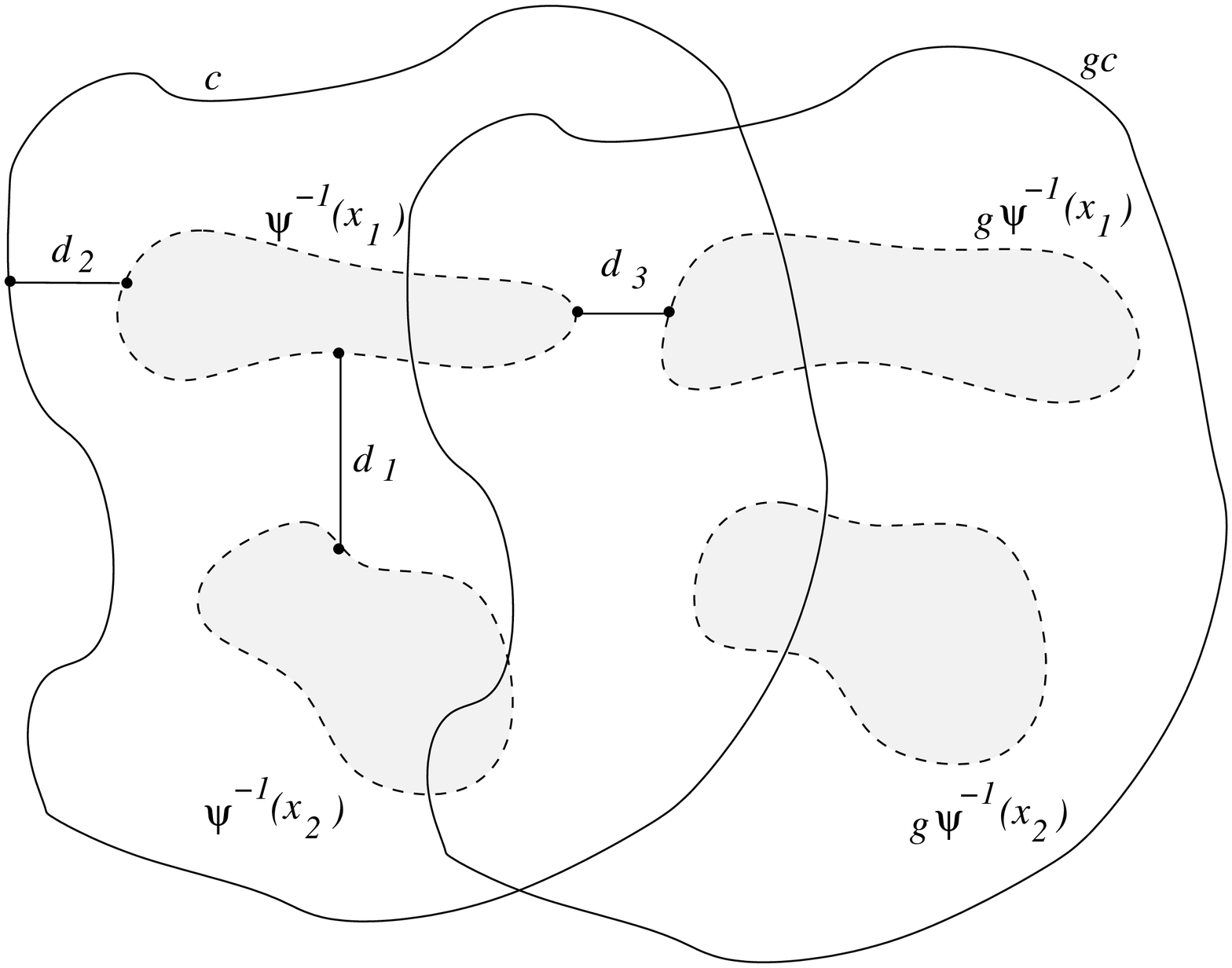}

Now, consider in $M$ an invariant metric and let $d_1$ be the distance between
$\psi^{-1}(x_1)$ and $\bigcup_{i\not=1} G\psi^{-1}(x_i)$. Let $d_2$ be the
distance between $\psi^{-1}(x_1)$ and the boundary of $c$ and $d_3$ the
distance between $\psi^{-1}(x_1)$ and $(G\psi^{-1}(x_1)-\psi^{-1}(x_1))\cap C$. Finally, take $d$ to be the minimum of
$\{d_1,d_2,d_3\}$. If we consider a finite triangulation of $c$ with mesh size less
than $d/3$, then no closed simplex of $c$  intersecting
$\psi^{-1}(x_1)$ intersects a simplex which touches $[\bigcup_{i}
(G\psi^{-1}(x_i) - \psi^{-1}(x_1)) \bigcup\partial{c}]\cap C = R_C$. Let $P_C$ be
the subpolyhedron of $c$ consisting of simplices which do not intersect $R_C$
and let $Q_C$ be the subpolyhedron of $P_C$ consisting of those simplices which do
not intersect $\psi^{-1}(x_1)$. Then $\sigma_{M}|_{Q_C}$ has no singularities and
again, by the same procedure used in Lemma 1.6 in \cite{brown}, we may extend it to $P_C$ with a finite number of
singularities, say, $\{y_1,y_2,...,y_m\}$. Since the metric on $M$ is invariant,
we may triangulate $g_{i}c$ in the same way we triangulate $c$ so that the
complexes $g_{i}P_C$ and $g_{i}Q_C$ will be the complexes corresponding to $P_C$ and
$Q_C$ for $g_{i}\psi^{-1}(x_1)$ in $g_{i} Q_C$. By doing so, the singularities of the extension
of  $\sigma_{M}|_{g_{i}Q_C}$ will be $\{g_{i}y_1, g_{i}y_2,..., g_{i}y_m\}$. Finally, we
extend this section to $M$ by making it agree with $\sigma_M$ outside
$G_{C}P_{C}$, where $G_C = \lbrace g\in G\,|\,gC=C \rbrace$. Repeating this procedure for $i>1$, we end up with an equivariant section extending $\sigma_A$ with a finite number of singular orbits,
$\{Gy_1,Gy_2,...,Gy_m\}$, lying in the closure of various components of $M-A$.

The next step is to reduce the set of singularities in such a way the closure of each
component of $M-A$ meet at most one singular orbit. For this, let $C$ be a component of $M-A$ and ${G_C y_1,...G_C y_r}$ be all
singularities in $C$. Consider $e$ a
closed $n$-cell in $C$ containing this entire set of singularities in its interior. Let
$e_1$ be another closed $n$-cell contained in the interior of $e$ such that
$y_1,y_2,...,y_r$ are in $e_1$ and $g_{l}e_1\bigcap g_{j}e_1=\emptyset$, for
$g_{l}$ and $g_{j}$ in $G_C$, $l\not=j$, as in the figure below.

\includegraphics[height=7cm]{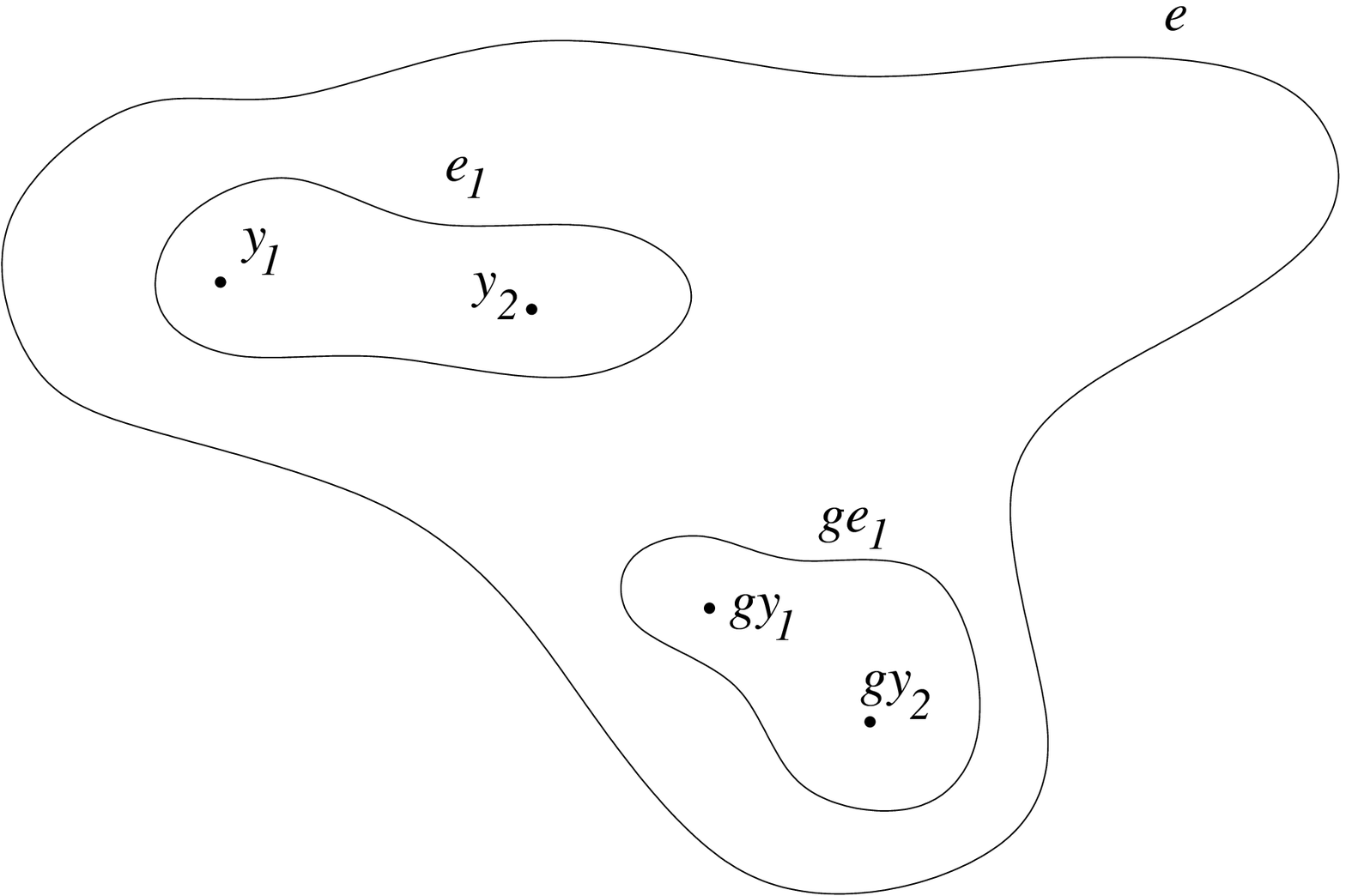}

Applying once more Lemma~\ref{lemma} for the $n$-cell $e_1$ we may reduce the set
$\{y_1,y_2,...,y_m\}$ to a single singularity, say $z$. Doing the same for the
cells $g_{j}e_1$, we end up with a cross section $\tau\colon M\lra T_M$ with
only $G_C z$ as singularities in $C$. Repeating this procedure for all other
components we are able to extend the path field to $M$ in such a way that the closure of each component $C$ of $M-A$ meet at most two singular orbits, one lying in the interior of $C$ and the other in its boundary. 

Now, for a component $C$ of $M-A$ let $G_{C}z\cup G_{C}w$ be its set of singularities, where $z\in C$ and $w\in (\overline{C})-C$. Then, it is possible to find $\vert G_{C} \vert$ cells touching $(\overline{C})-C$ only in $gw$, $g\in G_{C}$ and so that each of them contains only one pair of points of the form $gz$, $gw$, $g \in G_C$, as in the figure below.

\bg

\includegraphics[height=5.5cm]{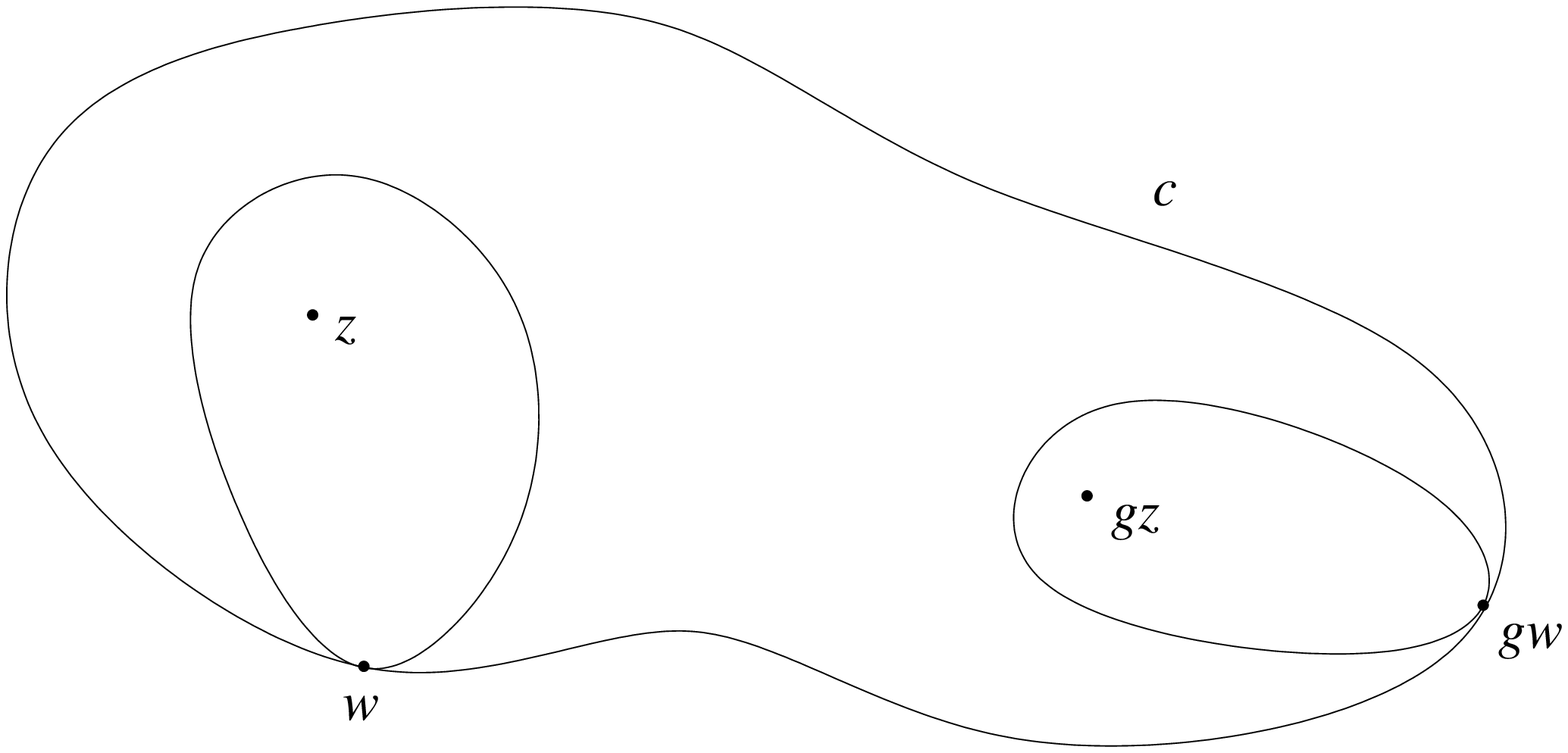}

Then applying Lemma~\ref{lemma} we see that we may change $\sigma$ in the interior of the cell containing $z$ and $w$ so that the only singularity left is $w$. Then repeating the procedure to the other cells, equivariantly, we complete the proof.  
\end{proof}

\begin{proposition}\label{induction-step}
Let $M$ be a $G-$manifold, $A\subset M^n$ an invariant submanifold so that $G$ acts freely on $M-A$. Assume that $A$ admits a $G$-path field without singularities. Then $M$ admits a $G$-path field with no singular orbits iff  $| \chi| (M-A)=0$.
\end{proposition}
\begin{proof}
Assume first $| \chi| (M-A)=0$ and let $\sigma$ be a $G$-path field on $M$ with a single singular orbit, say $Gx$, $x$ in a component $C$ of $M-A$. Take $f\colon M\lra M$, $f(x)=\sigma(x)(1)$. Then $f$ has only one fixed orbit in $C$, namely $G_{C}x$. Since $G$ is finite and $\chi(C)=0$ we have that the sum of the fixed point indices of $f$ at $gx$, $g$ in $G_C$, must vanish. Since the action of $G$ on $M$ is
locally smooth and the fixed points are isolated and lie in the same orbit, it
is not hard to see that they have all the same index, and therefore index zero. Because the action is free in $M-A$, we can find an Euclidean neighborhood, $U$, of $x$ in $C$ such that $gU\cap hU=\emptyset$, for all $g$ and $h$ in $G_C$. Applying exactly the same procedure as in the proof of Theorem 2.3 of \cite{brown}, we conclude that it is possible to construct a path field $\sigma'$ over $M$ so that it agrees with $\sigma$ in $M-U$ and has no singularities in $U$. Define $\tau\colon M\lra T_M$ to
agree with $\sigma$ in $M-\cup_{g\in G_C}gU$ and, for $y\in gU$,
$\tau(y)=g\sigma'(g^{-1}y)$. It is not difficult to see that $\tau$ is a
$G$-path field over $M$ without singular orbits in $C$. The proof is complete if
we repeat the same procedure for all other components of $M-A$.

Now, suppose $\sigma_A$ has no singularities and can be extended to $M$. Let
$\lbrace C_1,C_2,...,C_r \rbrace$ and $\lbrace A_1,A_2,...,A_l \rbrace$ be the
connected components of $M-A$ and $A$, respectively. Since, $\sigma|_{A_j}: A_j
\to A_j$ has no singularities, we have that $\chi(A_j) =0$, for all $j$.
Consider $D_i$ the union of the components of $A$ that meet the closure of
$C_i$. Set $\widetilde C_i = D_i\cup C_i$ and let $U_i$ be a tubular neighborhood of $D_i$ in
$\widetilde C_i$. Then $\chi(D_i) = 0$ and $\chi(U_i - D_i) = 0$, since $U_i-D_i$
fibers over $D_i$. Now, $\chi(C_i) = \chi(\widetilde C_i)- \chi(U_i) + \chi(U_i \cap
C_i) = \chi(\widetilde C_i)- \chi(D_i) +\chi(U_i-D_i) = \chi(\widetilde C_i)$.

Using the compactness of $\widetilde C_i$ and the fact that $\sigma(x)$ starts at
$x$, it is possible to find $t_i \in \lbrack 0,1 \rbrack$ so that
$\sigma(x)(\lbrack 0,t_i \rbrack)\subset \widetilde C_i$, for all $x$ in $\widetilde
C_i$. Therefore the map $f_i: \widetilde C_i \to \widetilde C_i$ given by $f_i (x)=
\sigma(x)(t_i)$ is fixed point free, and this implies that $\chi(\widetilde C_i) = 0$.
So we conclude that $\chi (C_i)=0$. This completes the proof.
\end{proof}

\begin{proposition}\label{one orbit type case}
Suppose $G$ acts on a manifold $M$ with only one orbit type $(H)$. Then it is possible to construct a $G$-path field on $M$ with a single singular orbit. Moreover, $M$
admits a $G$-path field with no singular orbits iff  $|\chi| (M_H)=0$.
\end{proposition}
\begin{proof}
Since $NH/H$ acts freely on $M^H$ then, by Proposition~\ref{relative free case}, it is possible to construct a $NH/H$-path field  
$\sigma'\colon M^H \lra M^H $ with only one $NH/H$- singular orbit. Since $G$ acts on $M$ with
only one orbit type $(H)$, we have that $M$ is $G$-homeomorphic to $G\times_{NH} M^H$, where $NH$ is the normalizer of $H$ (see \cite{bredon}). Define $\sigma\colon M\lra T_M$ by, $\sigma [g,x]=g\sigma'(x)$. It is not hard to see that $\sigma$ is an equivariant section with a  single singular orbit.

Now, assume $| \chi | (M_H)=0$. Then $| \chi | (M^H)=0$ and the section $\sigma'$ can be taken without singularities and so does $\sigma$. Finally, if $M$ admits a $G$-path field with no singular orbits, then $f\colon M\lra M$
given by $f(x)=\sigma(x)(1)$, has no fixed orbits. Therefore $f^H\colon M^H \lra
M^H $ has no fixed points, which implies that  $| \chi | (M^H)=0$ and the proof is done.
\end{proof}

\begin{theorem}\label{main-theorem}
Let $G$ be a finite group and $M$ a compact locally smooth $G$-space. Then  there exists a $G$-path field on $M$
having at most one singular orbit in the closure of each component of $M_H$. Moreover, $M$ admits a non singular $G$-path field iff $| \chi | (M_H)=0$, for all $H \le G$.
\end{theorem}
\begin{proof}
Consider $(H_1),(H_2),... (H_r)$ the orbit types of the $G$-action on $M$ ordered in a way that 
$(H_i)\subset (H_j)$ implies $j \leq i$. For each $i\in \{1,2,...r\}$, let $M_i = \{ x\in M \,|\, (G_x) = (H_j),j\leq i\}$. Then $M_1\subset M_2\subset ...\subset M_r$, $M_1 = M_{(H_1)}$, $M_r = M$ and 
$M_i - M_{i-1} = M_{(H_i)}$. Here, $M_{(H_j)} = \{x\in M | (G_x)= (H_j)\}$.

To prove the first part we will use induction on $r$. If $r=1$, then $M_1$ has only one orbit type,
namely, $(H_1)$. Therefore, Proposition~\ref{one orbit type case} implies
that $M_1$ admits a $G$-path field $\sigma_1$ with only one singular orbit.

Suppose we have succeeded  extending $\sigma_1$ to a $G$-path field,
$\sigma_{i-1}$, on $M_{i-1}$ so that the closure of each component of $M_{i-1} - M_{i-2} = M_{H_{i-1}}$ intersects at most one singular orbit. Take $N = M^{H_i} - (M_{i-1}\cap M^{H_i}) =M_{H_i}$. Since $NH_i/H_i$ acts freely on $N$, Proposition~\ref{one orbit type case} implies that we
are able to extend $\sigma_{i-1}|_{ M_{i-1} \cap M^{H_i}}$ to an $NH_i/H_i$-path
field, $\bar\sigma_i:M^{H_i} \to T_{M^{H_i}}$, without $NH_i/H_i$-singular
orbits.

Define $\widetilde\sigma_i: M^{(H_i)}\to T_{M^{(H_i)}}$ by $\widetilde\sigma_{i}(x) =l\bar\sigma_i(l^{-1}x)$, 
where $l\in G$ is such that $G_x \supset lH_{i}l^{-1}$.
Then, $\widetilde\sigma_i$ is a  well defined  $G$-path field
extending $\sigma_{i-1} |_{ M_{i-1} \cap M^{H_i}}$.

Now let $\sigma_{i}:M_{i} \to T_{M_i}$ coincide with
$\widetilde\sigma_i$, on $M^{(H_i)}$ and with $\sigma_{i-1}$ in $ M_{i-1}$. It is not
difficult to see that $\sigma_i$ is a $G$-path field extending
$\sigma_{i-1}$ with the desired property.

For the second part, assume $|\chi|(M_H)=0$, $\forall H\leq G$. The $G$-path field $\sigma$ constructed above can be taken without singular orbits by making use of Proposition~\ref{induction-step}, inductively on $\{(H)\}$.

Finally, if $M$ admits a non-singular $G$-path field $\sigma$ then
looking at  $M_{i-1} \cap M^{H_i}$ as an $NH_i$-submanifold of $M^{H_i}$ we
may repeat the proof  of Proposition~\ref{induction-step}, to obtain that $| \chi | (M_{H_i})=0$. The proof is complete.
\end{proof}

\section{$G$-Complete Invariance Property}

In this section, we study related problems concerning the fixed point theory for $G$-deformations. Recall that a $G$-space $X$ is said to have the $G$-CIP for $G$-deformations ($G$-CIPD), if for any nonempty closed invariant subset $A\subseteq X$, there exists a $G$-deformation $\lambda \sim_G 1_X$ such that $Fix\,\lambda =A$. In \cite{wong4}, necessary and sufficient conditions were given for smooth $G$-manifolds to possess the $G$-CIPD. As an application of Theorem~\ref{main-theorem}, we obtain the following

\begin{theorem}\label{G-CIP}
Let $G$ be a finite group and $M$ a compact locally smooth $G$-manifold. Suppose for each isotropy type $(H)$, $M^H$ has dimension at least $2$. Let $A\subset M$ be a non-empty closed invariant subset. Then the following are equivalent:
\begin{itemize}
\item There exists a $G$-deformation $\varphi:M\to M$ such that $A=Fix\,\varphi$.
\item $A\cap \bar C\ne \emptyset$ whenever $\chi(C)\ne 0$ for any connected component $C$ of $M_H$ and $\bar C$ denotes the closure of $C$ in $M^H$.
\end{itemize}
\end{theorem}
\begin{proof}
Suppose that there exists a $G$-deformation $\varphi:M\to M$ such that $Fix\,\varphi=A$. Let $C$ be a connected component of $M_H$ such that $\chi(C)\ne 0$. By excision, we have $H_*(M^H,M^H-C)=H_*(M_H,M_H-C)$ so that $\chi(M^H,M^H-C)=\chi(M_H,M_H-C)=\chi(C)\ne 0$. Using the relative Lefschetz fixed point theorem and Proposition 2.1 of \cite{wil}, we conclude that $\varphi$ must have a fixed point in the closure $\bar C$ of $C$ in $M^H$. Since $Fix \,\varphi =A$, it follows that $A\cap \bar C\ne \emptyset$.

Conversely, if $\chi(C)\ne 0$ then by Theorem~\ref{main-theorem} there exists a $G$-path field $\sigma$ such that $\sigma$ has one singular orbit in $G\bar C$. If $A\cap C\ne \emptyset$, this singular orbit lies in $A \cap GC$. If $A\cap C=\emptyset$, then the singular orbit must lie in $M^{>(H)}$ since $M_H$ is open and dense in $M^H$ and $M^H=M_H\cup (M^H\cap M^{>(H)})$, where $M^{>(H)}$ denotes the set of points of $M$ of isotropy type $(K)>(H)$. Now, let $\varphi(x)=\sigma(x)(t_x)$ where $t_x=d(x,A)$ and $d$ is a bounded $G$-invariant metric. Then, we have $Fix \,\varphi=A$ and $\varphi\sim 1_M$.
\end{proof}

\begin{remark} In the case where $M$ is a smooth $G$-manifold and $M^H/WH$ is connected for each $(H)$, the necessary and sufficient conditions obtained in \cite{wong4} can be derived from those of Theorem~\ref{G-CIP}. Our formulation resembles closely to case (B) of Theorem 1 of \cite{jiang2} except that Jiang considered connected components $C$ of $M^H$ instead. According to case (B) of \cite[Theorem 1]{jiang2}, there is a fixed point free $G$-deformation if $M^H$ is connected and $\chi(M^H)=0$ for all $(H)$. However, counter-examples have been found by D. Ferrario \cite{ferrario}. Therefore, Theorem~\ref{G-CIP} gives the correct necessary and sufficient conditions.
\end{remark} 

\begin{remark} The first example of a $G$-space $X$ in which each of the identity maps $1_{X^H}:X^H \to X^H$ is deformable to be fixed point free but $1_X$ is not $G$-deformable to be fixed point free was given by M. Izydorek and A. Vidal \cite{iz-vidal}. We would like to point out that one can easily modify their example (by taking the cartesion product with the unit interval) to give an example of a $G$-Wecken complex in the sense of \cite{wong3} such that the equivariant Euler characteristic used in \cite{wong3} is nonzero, that is, the identity is not equivariantly deformable to be fixed point free. The case for smooth $G$-manifolds was studied by D. Ferrario in \cite{ferrario} for more general $G$-maps.
\end{remark}

\end{document}